\documentclass{amsart}%
\usepackage{amsmath}
\usepackage{graphicx}
\usepackage{amsfonts}
\usepackage{amssymb}%
\newtheorem{theorem}{Theorem}
\theoremstyle{plain}

\newtheorem{definition}{Definition}
\newtheorem{example}{Example}

\numberwithin{equation}{section}
\begin{document}
\title{Higher Homotopy Hopf Algebras Found:\ A Ten Year Retrospective}
\author{Ronald Umble $^{1}$}
\address{Department of Mathematics\\
Millersville University of Pennsylvania\\
Millersville, PA. 17551}
\email{ron.umble@millersville.edu}
\thanks{$^{1}$ This research funded in part by a Millersville University faculty
research grant.}
\date{December 1, 2010}
\maketitle

\begin{abstract}
At the 1996 conference honoring Jim Stasheff in the year of his 60th birthday,
I initiated the search for $A_{\infty }$-bialgebras in a talk entitled 
\textquotedblleft In Search of Higher Homotopy Hopf Algebras.\textquotedblright\ 
\ The idea in that talk was to think of a DG
bialgebra as some (unknown) higher homotopy structure with trivial higher
order structure and apply a graded version of Gerstenhaber and Schack's
bialgebra deformation theory. \ Indeed, deformation cohomology, which detects some 
(but not all) $A_{\infty }$-bialgebra structure, motivated the definition 
given by S. Saneblidze and myself in 2004. 
\end{abstract}

\vspace{0.1in}
\begin{center}
\textit{To Murray Gerstenhaber and Jim Stasheff}
\end{center}

\section{Introduction}

In a preprint dated June 14, 2004, Samson Saneblidze and I announced the
definition of $A_{\infty}$-bialgebras \cite{SU3}, marking approximately six
years of collaboration that continues to this day. Unknown to us at the
time, $A_{\infty}$-bialgebras are ubiquitous and fundamentally important.
Indeed, over a field $F$, the bialgebra structure on the singular chains of
a loop space $\Omega X$ pulls back along a quasi-isomorphism $%
g:H_{\ast}\left( \Omega X;F\right) \rightarrow C_{\ast}\left( \Omega
X\right) $ to an $A_{\infty}$-bialgebra structure on homology that is unique
up to isomorphism \cite{SU5}.

Many have tried unsuccessfully to define $A_{\infty}$-bialgebras. The
illusive ingredient in the definition turned out to be an explicit diagonal $%
\Delta_{P}$ on the permutahedra $P=\sqcup_{n\geq1}P_{n}$, the first
construction of which was given by S. Saneblidze and myself in \cite{SU2}.
This paper is an account of the historical events leading up to the
discovery of $A_{\infty}$-bialgebras and the truly remarkable role played by 
$\Delta_{P}$ in this regard. Although the ideas and examples presented here
are quite simple, they represent and motivate general theory in \cite{SU2}, 
\cite{SU3}, \cite{SU4}, and \cite{SU5}.

Through their work in the theory of PROPs and the related area of infinity
Lie bialgebras, many authors have contributed indirectly to this work, most
notably M. Chas and D. Sullivan \cite{CS}, J-L. Loday 
\cite{Loday1}, M. Markl \cite{Mark1}, T. Pirashvili \cite{Pirashvili}, 
and B. Vallette \cite{Vallette}; for extensive
bibliographies see \cite{Sullivan} and \cite{Mark2}.

Several new results spin off of this discussion and are included here:\
Example \ref{ex1} in Section 3 introduces the first example of a bialgebra $%
H $ endowed with an $A_{\infty}$-algebra structure that is compatible with
the comultiplication. Example \ref{strict-A(4)} in Section 4, introduces the
first example of a \textquotedblleft non-operadic\textquotedblright\ $%
A_{\infty}$-bialgebra with a non-trivial operation $\omega^{2,2}:H^{%
\otimes2} \rightarrow H^{\otimes2}$. And in Section 5 we prove Theorem \ref%
{main}: \textit{Given a DG bialgebra} $\left( H,d,\mu,\Delta\right) $\textit{%
\ and a Gerstenhaber-Schack } $\mathit{2} $\textit{-cocycle }$\mu_{1}^{n}\in
Hom^{2-n}\left( H^{\otimes n},H\right) ,$\textit{\ }$n\geq3,$\textit{\ let} $%
H_{0}=\left( H\left[ \left[ t\right] \right] ,d,\mu,\Delta\right) .$ \textit{%
Then }$\left( H\left[ \left[ t\right] \right] ,d,\mu,\Delta,t\mu_{1}^{n}%
\right) $ \textit{is a linear deformation of }$H_{0}$ \textit{as a Hopf }
$A\left( n\right) $-\textit{algebra.}

\section{The Historical Context}

Two papers with far-reaching consequences in algebra and topology appeared
in 1963. In \cite{Gersten} Murray Gerstenhaber introduced the deformation
theory of associative algebras and in \cite{Stasheff1} Jim Stasheff
introduced the notion of an $A\left( n\right) $-algebra. Although the notion
of what we now call a \textquotedblleft non-$\Sigma$ operad%
\textquotedblright\ appears in both papers, this connection went unnoticed
until after Jim's visit to the University of Pennsylvania in 1983. Today,
Gerstenhaber's deformation theory and Stasheff's higher homotopy algebras
are fundamental tools in algebra, topology and physics. An extensive
bibliography of applications appears in \cite{MSS}.

By 1990, techniques from deformation theory and higher homotopy structures
had been applied by many authors, myself included \cite{Umble0}, \cite{LU},
to classify rational homotopy types with a fixed cohomology algebra. And it
seemed reasonable to expect that rational homotopy types with a fixed
Pontryagin algebra $H_{\ast}\left( \Omega X;\mathbb{Q}\right) $ could be
classified in a similar way. Presumably, such a theory would involve
deformations of DG bialgebras (DGBs) as some higher homotopy structure with
compatible $A_{\infty}$-algebra and $A_{\infty}$-coalgebra substructures,
but the notion of compatibility was not immediately clear and an appropriate
line of attack seemed illusive. But one thing was clear: If we apply a
graded version of Gerstenhaber and Schack's (G-S) deformation theory \cite%
{GS}, \cite{LM1}, \cite{LM2}, \cite{Umble1} and deform a DGB $H$ as some
(unknown) higher homotopy structure, new operations $\omega^{j,i}:H^{\otimes
i}\rightarrow H^{\otimes j}$ appear and their interactions with the deformed
bialgebra operations are partially detected by the differentials. While this
is but one small piece of a very large puzzle, it gave us a clue.

During the conference honoring Jim Stasheff in the year of his 60th
birthday, held at Vassar College in June 1996, I discussed this particular
clue in a talk entitled \textquotedblleft \textit{In Search of Higher
Homotopy Hopf Algebras}\textquotedblright\ (\cite{McCleary} p. xii).
Although G-S deformations of DGBs are less constrained than the $A_{\infty}$%
-bialgebras known today, they motivated the definition announced eight years
later.

Following the Vassar conference, forward progress halted. Questions of
structural compatibility seemed mysterious and inaccessible. Then in 1998,
Jim Stasheff ran across some related work by S. Saneblidze \cite{Saneblidze1}%
, of the A. Razmadze Mathematical Institute in Tbilisi, and suggested that I
get in touch with him. Thus began our long and fruitful collaboration. Over
the months that followed, Saneblidze applied techniques of homological
perturbation theory to solve the aforementioned classification problem \cite%
{Saneblidze2}, but the higher order structure in the limit is implicit and
the structure relations are inaccessible. In retrospect, this is not
surprising as explicit structure relations require explicit combinatorial
diagonals $\Delta_{P}$ on the permutahedra $P=\sqcup_{n\geq1}P_{n}$ and $%
\Delta_{K}$ on the associahedra $K=\sqcup_{n\geq2}K_{n}.$ But such diagonals
are difficult to construct and were unknown to us at the time. Indeed, one
defines the tensor product of $A_{\infty}$-algebras in terms of $\Delta_{K}$%
, and the search for a construction of $\Delta_{K}$ had remained a
long-standing problem in the theory of operads. We announced our
construction of $\Delta_{K}$ in 2000 \cite{SU1}; our construction of $%
\Delta_{P}$ followed a year or two later (see \cite{SU2}).

\section{Two Important Roles For $\Delta_{P}$}

The diagonal $\Delta_{P}$ plays two fundamentally important roles in the
theory of $A_{\infty}$-bialgebras. First, one builds the structure relations
from components of (co)free extensions of initial maps as higher
(co)derivations with respect to $\Delta_{P}$, and second, $\Delta_{P}$
specifies exactly which of these components to use.

To appreciate the first of these roles, recall the following definition
given by Stasheff in his seminal work on $A_{\infty }$-algebras in 1963 \cite%
{Stasheff1}: Let $A$ be a graded module, let $\{\mu ^{i}\in Hom^{i-2}\left(
A^{\otimes i},A\right) \}_{n\geq 1}$ be an arbitrary family of maps, and let 
$d$ be the cofree extension of $\Sigma \mu ^{i}$ as a coderivation of the
tensor coalgebra $T^{c}A$ (with a shift in dimension). Then $(A,\mu ^{i})$
is an \emph{$A_{\infty }$-algebra} if $d^{2}=0$; when this occurs, the
universal complex $\left( T^{c}A,d\right) $ is called the \emph{tilde-bar
construction} and the structure relations in $A$ are the homogeneous
components of $d^{2}=0$. Similarly, let $H$ be a graded module and let $%
\{\omega ^{j,i}\in Hom^{3-i-j}\left( H^{\otimes i},H^{\otimes j}\right)
\}_{i,j\geq 1}$, be an arbitrary family of maps. When $\left( H,\omega
^{j,i}\right) $ is an $A_{\infty }$-bialgebra, the map $\omega =\Sigma
\omega ^{j,i}$ uniquely extends to its \emph{biderivative }$d_{\omega }\in
End\left( TH\oplus T(H^{\otimes 2})\oplus \cdots \right) ,$ which is the sum
of various (co)free extensions of various subfamilies of $\left\{ \omega
^{j,i}\right\} $ as $\Delta _{P}$-(co)derivations (\cite{SU3}). And indeed,
the structure relations in $H$ are the homogeneous components of $d_{\omega
}^{2}=0$ with respect to an appropriate composition product.

To demonstrate the spirit of this, consider a free graded module $H$ of
finite type and an (arbitrary) map $\omega =\mu +\mu ^{3}+\Delta $ with
components $\mu :H^{\otimes 2}\rightarrow H,$ $\mu ^{3}:H^{\otimes
3}\rightarrow H,$ and $\Delta :H\rightarrow H^{\otimes 2}$. Extend $\Delta $
as a coalgebra map $\overline{\Delta }:T^{c}H\rightarrow T^{c}\left(
H^{\otimes 2}\right) ,$ extend $\mu +\mu ^{3}$ as a coderivation $%
d:T^{c}H\rightarrow T^{c}H$, and extend $\mu $ as an algebra map $\overline{%
\mu }:T^{a}\left( H^{\otimes 2}\right) \rightarrow T^{a}H.$ Finally, extend $%
\left( \mu \otimes 1\right) \mu $ and $\left( 1\otimes \mu \right) \mu $ as
algebra maps $f,g:T^{a}\left( H^{\otimes 3}\right) \rightarrow T^{a}H$, and
extend $\mu ^{3}$ as an $\left( f,g\right) $-derivation $\overline{\mu }%
^{3}:T^{a}\left( H^{\otimes 3}\right) \rightarrow T^{a}H.$ The components of
the biderivative in

\begin{center}
$d+\overline{\mu }+\overline{\mu }^{3}+\overline{\Delta }\in {\textstyle%
\bigoplus\limits_{p,q,r,s\geq 1}}Hom\left( \left( H^{\otimes p}\right)
^{\otimes q},\left( H^{\otimes r}\right) ^{\otimes s}\right) $
\end{center}

\noindent determine the structure relations. Let $\sigma _{r,s}:\left(
H^{\otimes r}\right) ^{\otimes s}\rightarrow \left( H^{\otimes s}\right)
^{\otimes r}$ denote the canonical permutation of tensor factors and define
a composition product $\circledcirc $ on homogeneous components $A$ and $B$
of $d+\overline{\mu }+\overline{\mu }^{3}+\overline{\Delta }$ by

\begin{center}
$A\circledcirc B=\left\{ 
\begin{array}{cl}
A\circ \sigma _{r,s}\circ B, & \text{if defined} \\ 
0, & \text{otherwise.}%
\end{array}%
\right. $
\end{center}

\noindent When $A\circledcirc B$ is defined, $\left( H^{\otimes r}\right)
^{\otimes s}$ is the target of $B$, and $\left( H^{\otimes s}\right)
^{\otimes r}$ is the source of $A.$ Then $\left( H,\mu ,\mu ^{3},\Delta
\right) $ is an $A_{\infty }$-infinity bialgebra if $d_{\omega }\circledcirc
d_{\omega }=0.$ Note that $\Delta \mu $ and $\left( \mu \otimes \mu \right)
\sigma _{2,2}\left( \Delta \otimes \Delta \right) $ are the homogeneous
components of $d_{\omega }\circledcirc d_{\omega }$ in $Hom\left( H^{\otimes
2},H^{\otimes 2}\right) ;$ consequently, $d_{\omega }\circledcirc d_{\omega
}=0$ implies the Hopf relation

\begin{center}
$\Delta\mu=\left( \mu\otimes\mu\right) \sigma_{2,2}\left( \Delta
\otimes\Delta\right) .$
\end{center}

\noindent Now if $\left( H,\mu,\Delta\right) $ is a bialgebra, the
operations $\mu_{t},$ $\mu_{t}^{3},$ and $\Delta_{t}$ in a G-S deformation
of $H$ satisfy 
\begin{equation}
\Delta_{t}\mu_{t}^{3}=\left[ \mu_{t}\left( \mu_{t}\otimes1\right)
\otimes\mu_{t}^{3}+\mu_{t}^{3}\otimes\mu_{t}\left( 1\otimes\mu_{t}\right) %
\right] \sigma_{2,3}\Delta_{t}^{\otimes3}  \label{five}
\end{equation}
and the homogeneous components of $d_{\omega}\circledcirc d_{\omega}=0$ in $%
Hom\left( H^{\otimes3},H^{\otimes2}\right) $ are exactly those in (\ref{five}%
). So this is encouraging.

Recall that the permutahedron $P_{1}$ is a point $0$ and $P_{2}$ is an
interval $01$. In these cases $\Delta_{P}$ agrees with the Alexander-Whitney
diagonal on the simplex:

\begin{center}
$\Delta_{P}\left( 0\right) =0\otimes0$ and $\Delta_{P}\left(
01\right) =0\otimes01+01\otimes1.$
\end{center}

\noindent If $X$ is an $n$-dimensional cellular complex, let $C_{\ast}\left(
X\right) $ denote the cellular chains of $X.$ When $X$ has a single top
dimensional cell, we denote it by $e^{n}.$ An $A_{\infty}$-algebra structure 
$\left\{ \mu^{n}\right\} _{n\geq2}$ on $H$ is encoded operadically by a
family of chain maps

\begin{center}
$\left\{ \xi:C_{\ast}\left( P_{n-1}\right) \rightarrow Hom\left( H^{\otimes
n},H\right) \right\} ,$
\end{center}

\noindent which factor through the map $\theta:C_{\ast}\left( P_{n-1}\right)
\rightarrow C_{\ast}\left( K_{n}\right) $ induced by cellular projection $%
P_{n-1}\rightarrow K_{n}$ given by A. Tonks \cite{Tonks} and satisfy $%
\xi\left( e^{n-2}\right) =\mu^{n}.$ The fact that

\begin{center}
$\left( \xi\otimes\xi\right) \Delta_{P}\left( e^{0}\right) =\mu\otimes\mu$ \
and

\vspace*{0.1in}$\left( \xi\otimes\xi\right) \Delta_{P}\left( e^{1}\right)
=\mu\left( \mu\otimes1\right) \otimes\mu^{3}+\mu^{3}\otimes\mu\left(
1\otimes\mu\right) $
\end{center}

\noindent are components of $\overline{\mu }$ and $\overline{\mu }^{3}$
suggests that we extend a given $\mu ^{n}$ as a higher derivation $\overline{%
\mu }^{n}:T^{a}\left( H^{\otimes n}\right) \rightarrow T^{a}H$ with respect
to $\Delta _{P}.$ Indeed, an $A_{\infty }$-bialgebra of the form $\left(
H,\Delta ,\mu ^{n}\right) _{n\geq 2}$ is defined in terms of the usual $%
A_{\infty }$-algebra relations together with the relations 
\begin{equation}
\Delta \mu ^{n}=\left[ \left( \xi \otimes \xi \right) \Delta _{P}\left(
e^{n-2}\right) \right] \sigma _{2,n}\Delta ^{\otimes n},  \label{six}
\end{equation}%
which define the compatibility of $\mu ^{n}$ and $\Delta $.

Structure relations in more general $A_{\infty }$-bialgebras of the form
\linebreak $\left( H,\Delta ^{m},\mu ^{n}\right) _{m,n\geq 2}$ are similar
in spirit and formulated in \cite{Umble2}. Special cases of the form $\left(
H,\Delta ,\Delta ^{n},\mu \right) $ with a single $\Delta ^{n}$ were studied
by H.J. Baues in the case $n=3$ \cite{Baues} and by A. Berciano and myself
with $n\geq 3$ \cite{Berciano}. Indeed, if $p$ is an odd prime and $n\geq 3$%
, these particular structures appear as tensor factors of the mod $p$
homology of an Eilenberg-Mac Lane space of type $K(\mathbb{Z},n)$.

Dually, $A_{\infty }$-bialgebras $\left( H,\Delta ,\mu ,\mu ^{n}\right) $
with a single $\mu ^{n}$ have a coassociative comultiplication $\Delta ,$ an
associative multiplication $\mu ,$ and $\xi \otimes \xi $ acts exclusively
on the primitive terms of $\Delta _{P}$ for lacunary reasons, in which case
relation (\ref{six}) reduces to 
\begin{equation}
\Delta \mu ^{n}=\left( f_{n}\otimes \mu ^{n}+\mu ^{n}\otimes f_{n}\right)
\sigma _{2,n}\Delta ^{\otimes n},  \label{seven}
\end{equation}%
where $f_{n}=\mu \left( \mu \otimes 1\right) \cdots \left( \mu \otimes
1^{\otimes n-2}\right) $. The first example of this particular structure now
follows.

\begin{example}
\label{ex1}Let $H$ be the primitively generated bialgebra $\Lambda\left(
x,y\right) $ with $\left\vert x\right\vert =1, \left\vert y\right\vert =2,$
and$\vspace*{0.1in}$\newline
$\hspace*{0.5in}\mu^{n}\left(
x^{i_{1}}y^{p_{1}}|\cdots|x^{i_{n}}y^{p_{n}}\right) =\left\{ 
\begin{array}{cl}
y^{p_{1}+\cdots+p_{n}+1}, & i_{1}\cdots i_{n}=1\ \text{and}\ p_{k}\geq1 \\ 
0, & \text{otherwise}.%
\end{array}
\right. \vspace*{0.1in}$\newline
One can easily check that $H$ is an $A_{\infty}$-algebra, and a
straightforward calculation together with the identity$\bigskip$\newline
$\left( 
\begin{array}{c}
p_{1}\!+\!\cdots\!+p_{n}\!+1 \\ 
i%
\end{array}
\right) =\!{\sum\limits_{s_{1}+\cdots+s_{n}=i-1}}\!\!\left( 
\begin{array}{c}
p_{1} \\ 
s_{1}%
\end{array}
\right) \!\cdots\!\left( 
\begin{array}{c}
p_{n} \\ 
s_{n}%
\end{array}
\right) +\!{\sum\limits_{s_{1}+\cdots+s_{n}=i}}\!\!\left( 
\begin{array}{c}
p_{1} \\ 
s_{1}%
\end{array}
\right) \!\cdots\!\left( 
\begin{array}{c}
p_{n} \\ 
s_{n}%
\end{array}
\right) \bigskip$\newline
verifies relation (\ref{seven}).
\end{example}

The second important role played by $\Delta_{P}$ is evident in $A_{\infty}$%
-bialgebras in which $\omega^{n,m}$ is non-trivial for some $m,n>1.$ Just as
an $A_{\infty}$-algebra structure on $H$ is encoded operadically, an $%
A_{\infty}$-bialgebra structure on $H$ is encoded \emph{matradically }by a
family of chain maps

\begin{center}
$\left\{ \varepsilon:C_{\ast}\left( KK_{n,m}\right) \rightarrow
Hom^{3-m-n}\left( H^{\otimes m},H^{\otimes n}\right) \right\} $
\end{center}

\noindent over contractible polytopes $KK=\sqcup _{m,n\geq 1}KK_{n,m}$,
called \emph{biassociahedra}, with single top dimensional cells $e^{m+n-3}$ such
that $\varepsilon \left( e^{m+n-3}\right) =\omega ^{n,m}.$ Note that $%
KK_{n,1}=KK_{1,n}$ is the associahedron $K_{n}$ \cite{SU4}. Let $%
M=\{M_{n,m}= $ $Hom(H^{\otimes m},H^{\otimes n})\}$ and let $\Theta =\left\{
\theta _{m}^{n}=\omega ^{n,m}\right\} .$ The $A_{\infty }$\emph{-bialgebra
matrad} $\mathcal{H}_{\infty }$ is realized by $C_{\ast }\left( KK\right) $
and is a proper submodule of the free PROP $M$ generated by $\Theta .$ The 
\emph{matrad product }$\gamma $ on $\mathcal{H}_{\infty }$ is defined in
terms of $\Delta _{P},$ and a monomial $\alpha $ in the free PROP $M$ is a
component of a structure relation if and only if $\alpha \in \mathcal{H}%
_{\infty }$.

More precisely, in \cite{Mark1} M. Markl defined the submodule $S$ of \emph{%
special elements} in PROP $M$ whose additive generators are monomials $%
\alpha $ expressed as \textquotedblleft elementary
fractions\textquotedblright 
\begin{equation}
\alpha =\frac{\alpha _{p}^{y_{1}}\cdots \alpha _{p}^{y_{q}}}{\alpha
_{x_{1}}^{q}\cdots \alpha _{x_{p}}^{q}}  \label{fraction}
\end{equation}%
in which $\alpha _{x_{i}}^{q}$ and $\alpha _{p}^{y_{j}}$ are additive
generators of $S$ and the $j^{th}$ output of $\alpha _{x_{i}}^{q}$ is linked
to the $i^{th}$ input of $\alpha _{p}^{y_{j}}$ (here juxtaposition denotes
tensor product). Representing $\theta _{m}^{n}$ graphically as a double
corolla (see Figure 1), a general decomposable $\alpha $ is represented by a
connected non-planar graph in which the generators appear in order from
left-to-right (see Figure 2). The matrad $\mathcal{H}_{\infty }$ is a proper
submodule of $S$ and the matrad product $\gamma $ agrees with the
restriction of Markl's fraction product to $\mathcal{H}_{\infty }$. \vspace*{%
0.1in}

\begin{center}
$%
\begin{array}{c}
\theta_{m}^{n}=%
\end{array}%
\begin{array}{c}
n \\ 
{\includegraphics[
height=0.482in,
width=0.4915in
]{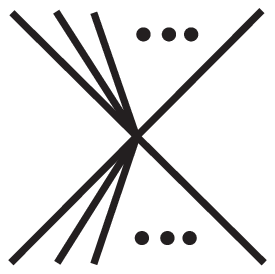}} \\ 
m%
\end{array}
$\vspace*{0.1in}

Figure 1.\vspace*{0.1in}
\end{center}

The diagonal $\Delta_{P}$ acts as a filter and admits certain elementary
fractions as additive generators of $\mathcal{H}_{\infty}$. In dimensions $0$
and $1,$ the diagonal $\Delta_{P}$ is expressed graphically in terms of
up-rooted planar rooted trees (with levels) by

\begin{center}
$\Delta_{P}(\ 
\raisebox{-0.0363in}{\includegraphics[
height=0.1773in,
width=0.1816in
]{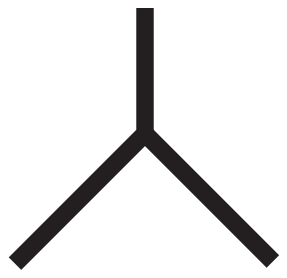}} \ )=\ 
\raisebox{-0.0363in}{\includegraphics[
height=0.1773in,
width=0.1816in
]{T2.eps}} \otimes\ 
\raisebox{-0.0363in}{\includegraphics[
height=0.1773in,
width=0.1816in
]{T2.eps}} $ \ \ and \ $\Delta_{P}(\ 
\raisebox{-0.0571in}{\includegraphics[
height=0.2101in,
width=0.2145in
]{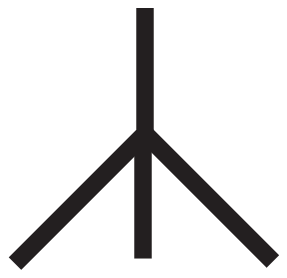}} \ )=\ 
\raisebox{-0.0623in}{\includegraphics[
height=0.2101in,
width=0.2145in
]{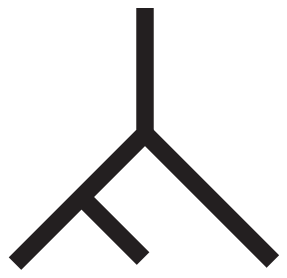}} \ \otimes\ 
\raisebox{-0.0467in}{\includegraphics[
height=0.2101in,
width=0.2145in
]{T3.eps}} \ +\ 
\raisebox{-0.0467in}{\includegraphics[
height=0.2101in,
width=0.2145in
]{T3.eps}} \ \otimes\ 
\raisebox{-0.0519in}{\includegraphics[
height=0.2101in,
width=0.2145in
]{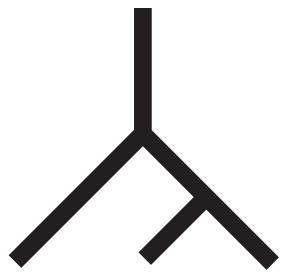}} \ .$
\end{center}

\noindent Define $\Delta_{P}^{(0)}={1};$ for each $k\geq1,$ define $\Delta
_{P}^{(k)}=\left( \Delta_{P}\otimes{1}^{\otimes k-1}\right) \Delta
_{P}^{(k-1)}$ and view each component of $\Delta_{P}^{(k)}(\theta_{q}^{1})$
as a $\left( q-2\right) $-dimensional subcomplex of $\left( P_{q-1}\right)
^{\times k+1},$ and similarly for $\Delta_{P}^{(k)}\left(
\theta_{1}^{q}\right) $.

The elements $\theta_{1}^{1},$ $\theta_{2}^{1},$ and $\theta_{1}^{2}$
generate two elementary fractions in $M_{2,2}$ each of dimension zero,
namely,

\begin{center}
$\alpha_{2}^{2}\mathtt{\ }=\mathtt{\ } 
\raisebox{-0.1773in}{\includegraphics[
height=0.416in,
width=0.1903in
]{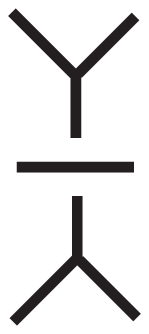}} $ \ and \ $\alpha_{11}^{11}\mathtt{\ }=\mathtt{\ 
\raisebox{-0.2084in}{\includegraphics[
height=0.4705in,
width=0.4869in
]{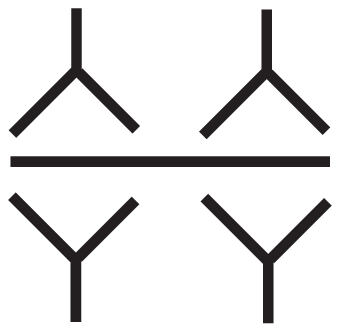}} \ .}$
\end{center}

\noindent Define $\partial\left( \theta_{2}^{2}\right)
=\alpha_{2}^{2}+\alpha_{11}^{11},$ and label the edge and vertices of the
interval $KK_{2,2}$ by $\theta_{2}^{2},$ $\alpha_{2}^{2}$ and $%
\alpha_{11}^{11},$ respectively. Continuing inductively, the elements $%
\theta_{1}^{1},$ $\theta_{2}^{1},$ $\theta_{1}^{2},$ $\theta_{2}^{2},$ $%
\alpha_{2}^{2},$ and $\alpha_{11}^{11}$ generate 18 fractions in $M_{2,3}$
-- one in dimension 2, nine in dimension 1 and eight in dimension 0. Of
these, 14 label the edges and vertices of the heptagon $KK_{2,3}.$ Since the
generator $\theta_{3}^{2}$ must label the 2-face, we wish to discard the
2-dimensional decomposable

\begin{center}
$e\mathtt{\ }=\mathtt{\ 
\raisebox{-0.1513in}{\includegraphics[
height=0.4246in,
width=0.5716in
]{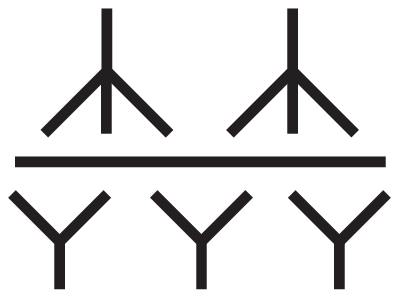}} }$
\end{center}

\noindent and the appropriate components of its boundary. Note that $e$ is a
square whose boundary is the union of four edges%
\begin{equation}
{\includegraphics[
height=0.4281in,
width=3.1816in
]{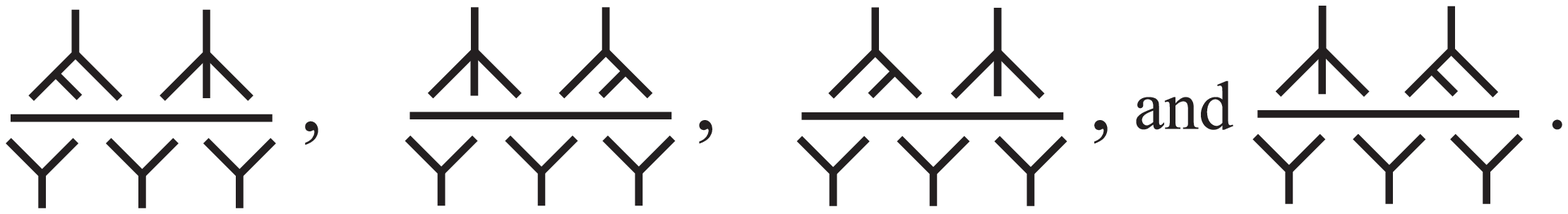}}  \label{four edges}
\end{equation}
Of the five fractions pictured above, only the first two in (\ref{four edges}%
) have numerators and denominators that are components of $%
\Delta_{P}^{\left( k\right) }\left( P\right) $ (numerators are components of 
$\Delta _{P}^{\left( 1\right) }(\theta_{3}^{1})$ and denominators are
exactly $\Delta_{P}^{\left( 2\right) }(\theta_{1}^{2})$). Our selection rule
admits only these two particular fractions, leaving seven 1-dimensional
generators to label the edges of $KK_{2,3}$ (see Figure 2). Now linearly
extend the boundary map $\partial$ to the seven admissible 1-dimensional
generators and compute the seven 0-dimensional generators labeling the
vertices of $KK_{2,3}.$ Since the 0-dimensional generator

\begin{center}
${\includegraphics[
height=0.4246in,
width=0.5716in
]{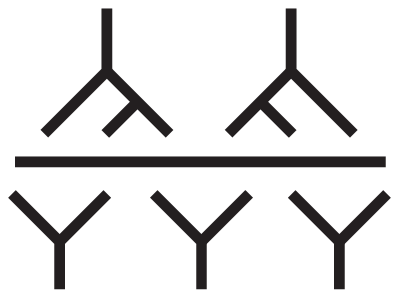}} $
\end{center}

\noindent is not among them, we discard it. \vspace{.2in}

\begin{center}
{\includegraphics[
height=2.444in,
width=2.5806in
]{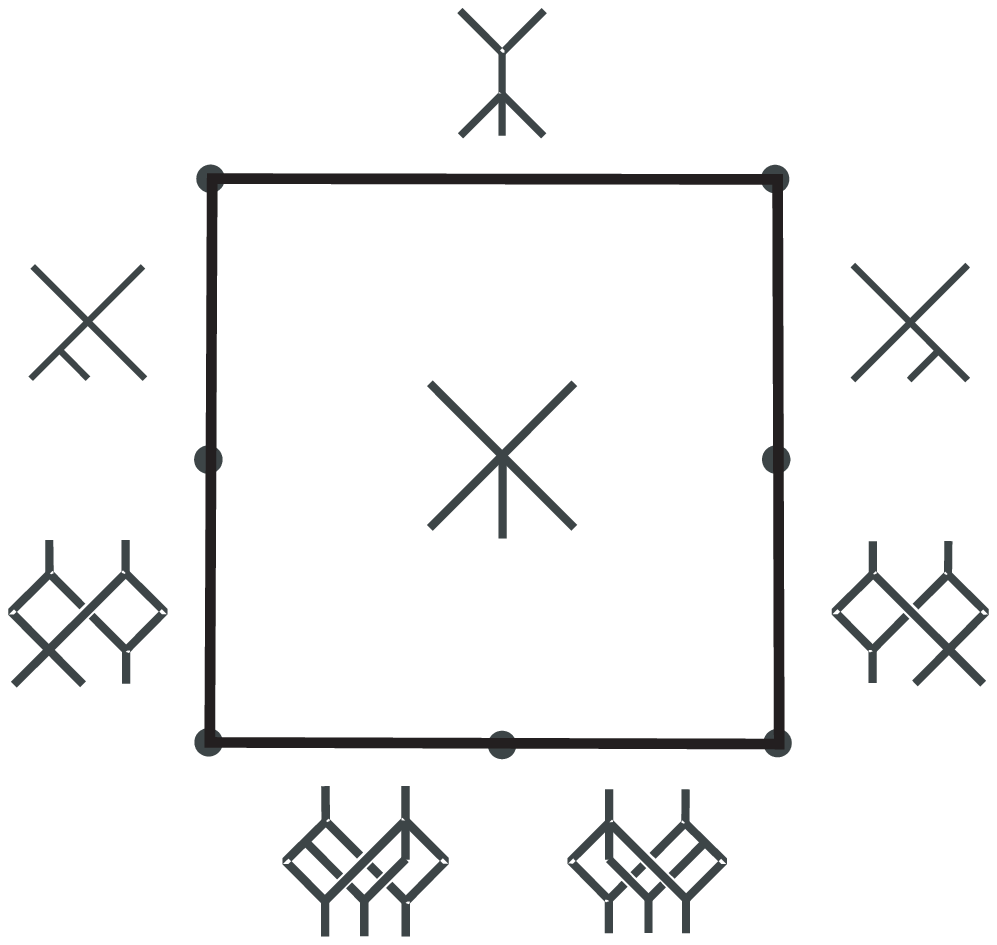}} \vspace*{0.1in}

Figure 2:\ The biassociahedron $KK_{2,3}.\vspace*{0.1in}$
\end{center}

Subtleties notwithstanding, this process continues indefinitely and produces 
$\mathcal{H}_{\infty}$ as an explicit free resolution of the bialgebra
matrad $\mathcal{H=}\left\langle
\theta_{1}^{1},\theta_{2}^{1},\theta_{1}^{2}\right\rangle $ in the category
of matrads. We note that in \cite{Mark1}, M. Markl makes arbitrary choices
(independent of our selection rule) to construct the polytopes $%
B_{m}^{n}=KK_{n,m}$ for $m+n\leq6$. In this range, it is enough to consider
components of the diagonal $\Delta_{K}$ on the associahedra.

We conclude this section with a brief review of our diagonals $\Delta_{P}$
and $\Delta_{K}$ (up to sign); for details see \cite{SU2}. Alternative
constructions of $\Delta_{K}$ were subsequently given by Markl and Shnider 
\cite{MS} and J.L. Loday \cite{Loday2} (in this volume). Let $\underline {n}%
=\{1,2,\dots,n\},$ $n\geq1.$ A matrix $E$ with entries from $\left\{
0\right\} \cup\underline{n}$ is a \emph{step matrix} if:

\begin{itemize}
\item Each element of $\underline{n}$ appears as an entry of $E$ exactly
once.

\item The elements of $\underline{n}$ in each row and column of $E$ form an
increasing contiguous block.

\item Each diagonal parallel to the main diagonal of $E$ contains exactly
one element of $\underline{n}$.
\end{itemize}

\noindent Right-shift and down-shift matrix transformations, which include
the identity (a trivial shift), act on step matrices and produce \emph{%
derived matrices.}

Given a $q\times p$ integer matrix $M=\left( m_{ij}\right) ,$ choose proper
subsets $S_{i}\subset \left\{ \text{non-zero}\right. $ $\left.\text{entries in row } i\right\}$ and $T_{j}\subset \left\{\text{non-zero entries in column } j\right\}$, and
define \emph{down-shift} and \emph{right-shift} operations\emph{\ }$%
D_{S_{i}} $ and $R_{T_{j}}$ on $M$ as follows:
\begin{itemize}
\item If $S_{i}\neq \varnothing ,$ $\max$ row$(i+1)<\min S_{i}=m_{ij},$ and $m_{i+1,k}=0$ for all $k\geq j$, then $%
D_{S_{i}}M$ is the matrix obtained from $M$ by interchanging each $m_{ik}\in
S_{i}$ with $m_{i+1,k};$ otherwise $D_{S_{i}}M=M.$
\item If $T_{j}\neq \varnothing ,$ $\max$ col$(j+1)<\min T_{j}=m_{ij},$ and $m_{k,j+1}=0$ for all $k\geq i,$ then\emph{\ }$%
R_{T_{j}}M$ is the matrix obtained from $M$ by interchanging each $%
m_{k,j}\in T_{j}$ with $m_{k,j+1};$ otherwise $R_{Tj}M=M.$
\end{itemize}
\noindent Given a $q\times p$ step matrix $E\ $together with subsets $%
S_{1},\ldots ,S_{q}$ and $T_{1},\ldots ,T_{p}$ as above, there is the \emph{%
derived matrix}%
\begin{equation*}
R_{T_{p}}\cdots R_{T_{2}}R_{T_{1}}D_{S_{q}}\cdots D_{S_{2}}D_{S_{1}}E.
\end{equation*}%
In particular, step matrices are derived matrices under the trivial action
with $S_{i}=T_{j}=\varnothing $ for all $i,j$.

Let $a=A_{1}|A_{2}|\cdots|A_{p}$ and $b=B_{q}|B_{q-1}|%
\cdots|B_{1}$ be partitions of $\underline{n}$. The pair $a\times b$ is an 
\emph{$(p,q)$-complementary pair} (CP) if $B_{i}$ and $A_{j}$ are the rows
and columns of a $q\times p$ derived matrix. Since faces of $P_{n}$ are
indexed by partitions of $\underline{n},$ and CPs are in one-to-one
correspondence with derived matrices, each CP is identified with some
product face of $P_{n}\times P_{n}.$

\begin{definition}
Define $\Delta_{P}(e^{0})=e^{0}\otimes e^{0}$. Inductively, having defined $%
\Delta_{P}$ on $C_{\ast}(P_{k+1})$ for all $0\leq k\leq n-1$, define $%
\Delta_{P}$ on $C_{n}(P_{n+1})$ by 
\begin{equation*}
\Delta_{P}(e^{n})=\sum_{\substack{ (p,q)\text{-CPs }u\times v  \\ p+q=n+2}}%
\pm\ u\otimes v,
\end{equation*}
and extend multiplicatively to all of $C_{\ast}(P_{n+1})$.
\end{definition}

The diagonal $\Delta_{P}$ induces a diagonal $\Delta_{K}$ on $C_{\ast}\left(
K\right) $. Recall that faces of $P_{n}$ in codimension $k$ are indexed by
planar rooted trees with $n+1$ leaves and $k+1$ levels (PLTs), and
forgetting levels defines the cellular projection $\theta:P_{n}\rightarrow
K_{n+1}$ given by A. Tonks \cite{Tonks}. Thus faces of $P_{n}$ indexed by
PLTs with multiple nodes in the same level degenerate under $\theta$, and
corresponding generators lie in the kernel of the induced map $%
\theta:C_{\ast}\left( P_{n}\right) \rightarrow C_{\ast}\left( K_{n+1}\right)$%
. The diagonal $\Delta_{K}$ is given by $\Delta_{K}\theta=(\theta\otimes%
\theta)\Delta_{P}.$

\section{Deformations of DG Bialgebras as $A\left( n\right) $-Bialgebras}

The discussion above provides the context to appreciate the extent to which
G-S deformation theory motivates the notion of an of $A_{\infty}$-bialgebra.
We describe this motivation in this section. In retrospect, the bi(co)module
structure encoded in the G-S differentials controls some (but not all) of
the $A_{\infty}$-bialgebra structure relations. For example, all structure
relation in $A_{\infty}$-bialgebras of the form $\left(
H,d,\mu,\Delta,\mu^{n}\right) $ are controlled except 
\begin{equation}
{\displaystyle\sum\limits_{i=0}^{n-1}}\left( -1\right) ^{i\left( n+1\right)
}\mu^{n}\left( 1^{\otimes i}\otimes\mu^{n}\otimes1^{\otimes n-i-1}\right) =0,
\label{relation}
\end{equation}
which measures the interaction of $\mu^{n}$ with itself. Nevertheless, such
structures admit an $A\left( n\right) $-algebra substructure and their
single higher order operation $\mu^{n}$ is compatible with $\Delta.$ Thus we
refer to such structures here as \emph{Hopf }$A\left( n\right) $\emph{%
-algebras. }General G-S deformations of DGBs, referred to here as \emph{%
quasi-}$A\left( n\right) $\emph{-bialgebras, }are \textquotedblleft
partial\textquotedblright\ $A\left( n\right) $-bialgebras in the sense that
all structure relations involving multiple higher order operations are out
of control.

\subsection{$A\left( n\right) $-Algebras and Their Duals}

The signs in the following definition were given in \cite{SU2} and differ
from those given by Stasheff in \cite{Stasheff1}. We note that either choice
of signs induces an oriented combinatorial structure on the associahedra,
and these structures are are equivalent. Let $n\in \mathbb{N}\cup\{\infty\}.$

\begin{definition}
\label{A(4)-defn} An $A\left( n\right) $\textbf{-algebra} is a graded module $A
$ together with structure maps $\{\mu^{k}\in Hom^{2-k}\left( A^{\otimes
k},A\right) \}_{1\leq k < n+1}$ that satisfy the relations

${\displaystyle\sum\limits_{j=0}^{k-1}}\ {\displaystyle\sum%
\limits_{i=0}^{k-j-1}}\left( -1\right) ^{j\left( i+1\right) }\mu^{k-j}\left(
1^{\otimes i}\otimes\mu^{j+1}\otimes1^{\otimes k-j-1-i}\right) =0$

\noindent for each $k<n+1$. Dually, an $A\left( n\right) $-{coalgebra}
is a graded module $C$ together with structure maps $\{\Delta^{k}\in
Hom^{2-k}(C,C^{\otimes j})\}_{1\leq k < n+1}$ that satisfy the relations

${\displaystyle\sum\limits_{j=0}^{k-1}}\ {\displaystyle\sum%
\limits_{i=0}^{k-j-1}}\left( -1\right) ^{j\left( k+i+1\right) }\left(
1^{\otimes i}\otimes\Delta^{j+1}\otimes1^{\otimes k-j-1-i}\right)
\Delta^{k-j}=0$

\noindent for each $k<n+1$.
\end{definition}

An $A\left( n\right) $-algebra is \emph{strict }if $n < \infty$ and $%
\mu^{n}=0.$ A \emph{simple }$A\left( n\right) $-algebra is a strict $A\left(
n+1\right) $-algebra of the form $\left( A,d,\mu,\mu^{n}\right) $; in
particular, a simple\emph{\ }$A\left( 3\right) $-algebra is a strict $%
A\left( 4\right) $-algebra in which

\begin{enumerate}
\item[i.] $d$ is both a differential and a derivation of $\mu,$

\item[ii.] $\mu$ is homotopy associative and $\mu^{3}$ is an associating
homotopy:\vspace*{0.1in}\newline
$\hspace*{0.3in}d\mu^{3}+\mu^{3}\left( d\otimes1\otimes1+1\otimes
d\otimes1+1\otimes1\otimes d\right) =\mu\left( \mu\otimes1\right) -\mu\left(
1\otimes\mu\right) ,$\vspace*{0.1in}

\item[iii.] $\mu$ and $\mu^{3}$ satisfy a strict pentagon condition:\vspace *%
{0.1in}\newline
$\hspace*{0.3in}\mu^{3}\left( \mu\otimes1\otimes1-1\otimes
\mu\otimes1+1\otimes1\otimes\mu\right) =\mu\left( 1\otimes\mu^{3}+\mu
^{3}\otimes1\right) .$\vspace*{0.1in}
\end{enumerate}

\subsection{Deformations of DG Bialgebras}

In \cite{GS}, M. Gerstenhaber and S. D. Schack defined the cohomology of an 
\emph{ungraded} bialgebra by joining the dual cohomology theories of G.
Hochschild \cite{Hochschild} and P. Cartier \cite{Cartier}. This
construction was given independently by A. Lazarev and M. Movshev in \cite%
{LM1}. The \emph{G-S cohomology of} $H$ reviewed here is a straight-forward
extension to the graded case and was constructed in \cite{LM2} and \cite%
{Umble1}.

Let $\left( H,d,\mu,\Delta\right) $ be a connected DGB. We assume $%
\left\vert d\right\vert =1,$ although one could assume $\left\vert
d\right\vert =-1$ equally well. For detailed derivations of the formulas
that follow see \cite{Umble1}. For each $i\geq1,$ the $i$\textit{-}\emph{%
fold bicomodule tensor power of}\textit{\ }$H$\textit{\ }is the $H$%
-bicomodule $H^{\underline{\otimes}i}=(H^{\otimes i},\lambda_{i},\rho_{i})$
with left and right coactions given by

\begin{center}
$\lambda_{i}=\left[ \mu\left( \mu\otimes1\right) \cdots(\mu\otimes
1^{\otimes i-2})\otimes1^{\otimes i}\right] \sigma_{2,i}\Delta^{\otimes i}\ $
and\vspace*{0.1in}

$\rho_{i}=\left[ 1^{\otimes i}\otimes\mu\left( 1\otimes\mu\right)
\cdots(1^{\otimes i-2}\otimes\mu)\right] \sigma_{2,i}\Delta^{\otimes i}$.
\end{center}

\noindent When $f:H^{\otimes i}\rightarrow H^{\otimes\ast},$ there is the
composition

\begin{center}
$\left( 1\otimes f\right) \lambda_{i}=\left[ \mu\left( \mu\otimes1\right)
\cdots(\mu\otimes1^{\otimes i-2})\otimes f\right] \sigma_{2,i}\Delta^{%
\otimes i}$.
\end{center}

\noindent Dually, for each $j\geq1,$ the $j$\emph{-fold bimodule tensor
power of}\textit{\ }$H$\textit{\ }is the $H$-bimodule $H^{\overline{\otimes}%
j}=(H^{\otimes j},\lambda^{j},\rho^{j})$ with left and right actions given by

\begin{center}
$\lambda^{j}=\mu^{\otimes j}\sigma_{j,2}\left[ (\Delta\otimes1^{\otimes
j-2})\cdots\left( \Delta\otimes1\right) \Delta\otimes1^{\otimes j}\right] \ $
and\vspace*{0.1in}

$\rho^{j}=\mu^{\otimes j}\sigma_{j,2}\left[ 1^{\otimes j}\otimes(1^{\otimes
j-2}\otimes\Delta)\cdots\left( 1\otimes\Delta\right) \Delta\right] $ .
\end{center}

\noindent When $g:H^{\otimes\ast}\rightarrow H^{\otimes j},$ there is the
composition

\begin{center}
$\lambda^{j}\left( 1\otimes g\right) =\mu^{\otimes j}\sigma_{j,2}\left[
(\Delta\otimes1^{\otimes j-2})\cdots\left( \Delta\otimes1\right)
\Delta\otimes g\right] $ .
\end{center}

Let $\mathbf{k}$ be a field. Extend $d,$ $\mu$ and $\Delta$ to $\mathbf{k}%
\left[ \left[ t\right] \right] $-linear maps and obtain a $\mathbf{k}\left[ %
\left[ t\right] \right] $-DGB $H_{0}=(H\left[ \left[ t\right] \right]
,d,\mu,\Delta).$ We wish to deform $H_{0}$ as an $A\left( n\right) $%
-structure of the form

\begin{center}
$H_{t}=\left( H\left[ \left[ t\right] \right] ,d_{t}=\omega_{t}^{1,1},%
\mu_{t}=\omega_{t}^{1,2},\Delta_{t}=\omega_{t}^{2,1},\omega_{t}^{j,i}\right)
_{i+j=n+1}$ \ where

\vspace*{0.1in}$\omega_{t}^{j,i}=\sum\limits_{k=0}^{\infty}t^{k}%
\omega_{k}^{j,i}\in Hom^{3-i-j}\left( H^{\underline{\otimes}i},H^{\overline{%
\otimes}j}\right) $,

\vspace*{0.1in}$\omega_{0}^{1,1}=d,$ $\omega_{0}^{1,2}=\mu,$ $%
\omega_{0}^{2,1}=\Delta,\mathtt{\ and}$ $\omega_{0}^{j,i}=0$.
\end{center}

\noindent Deformations of $H_{0}$ are controlled by the \emph{G-S }$n$-\emph{%
complex,} which we now review. For $k\geq1,$ let

\begin{center}
$%
\begin{array}{cll}
d_{(k)} & = & \sum\limits_{i=0}^{k-1}1^{\otimes i}\otimes d\otimes1^{\otimes
k-i-1} \\ 
&  &  \\ 
\partial_{(k)} & = & \sum\limits_{i=0}^{k-1}(-1)^{i}1^{\otimes i}\otimes
\mu\otimes1^{\otimes k-i-1} \\ 
&  &  \\ 
\delta_{(k)} & = & \sum\limits_{i=0}^{k-1}(-1)^{i}1^{\otimes i}\otimes
\Delta\otimes1^{\otimes k-i-1}.%
\end{array}
$
\end{center}

\noindent These differentials induce strictly commuting differentials $d,$ $%
\partial,$ and $\delta$ on the trigraded module $\{Hom^{p}(H^{\underline{%
\otimes}i},H^{\overline{\otimes}j})\},$ which act on an element $f$ in
tridegree $\left( p,i,j\right) $ by

\begin{center}
$%
\begin{array}{lll}
d(f) & = & d_{(j)}f-(-1)^{p}fd_{(i)} \\ 
&  &  \\ 
\partial(f) & = & \lambda^{j}(1\otimes f)-f\partial_{(i)}-(-1)^{i}\rho
^{j}(f\otimes1) \\ 
&  &  \\ 
\delta(f) & = & (1\otimes f)\lambda_{i}-\delta_{(j)}f-(-1)^{j}(f\otimes
1)\rho_{i}.%
\end{array}
$
\end{center}

\noindent The submodule of \emph{total G-S }$r$\emph{-cochains on }$H$ is

\begin{center}
$C_{GS}^{r}(H,H)=\bigoplus\limits_{p+i+j=r+1}Hom^{p}(H^{\underline{\otimes}%
i},H^{\overline{\otimes}j})$
\end{center}

\noindent and the total differential $D$ on a cochain $f$ in tridegree $%
\left( p,i,j\right) $ is given by

\begin{center}
$D\left( f\right) =\left[ \left( -1\right) ^{i+j}d+\partial+\left( -1\right)
^{i}\delta\right] \left( f\right) $ ,
\end{center}

\noindent where the sign coefficients are chosen so that (1) $D^{2}=0,$ (2)
structure relations (ii) and (iii) in Definition \ref{A(4)-defn} hold and
(3) the restriction of $D$ to the submodule of $r$-cochains in degree $p=0$
agrees with the total (ungraded) G-S differential. The \emph{G-S cohomology
of }$H$\emph{\ with coefficients in} $H$ is given by

\begin{center}
$H_{GS}^{\ast}\left( H,H\right) =H_{\ast}\left\{ C_{GS}^{r}\left( H,H\right)
,D\right\} .$
\end{center}

Identify $Hom^{p}(H^{\underline{\otimes}i},H^{\overline{\otimes}j})$ with
the point $\left( p,i,j\right) $ in $\mathbb{R}^{3}.$ Then the \emph{G-S }$n$%
-\emph{complex} is that portion of the G-S complex in the region $x\geq2-n$
and the submodule of total $r$-cochains in the $n$-complex is

\begin{center}
$C_{GS}^{r}(H,H;n)=\bigoplus\limits_{p=r-i-j+1\geq2-n}Hom^{p}(H^{\underline {%
\otimes}i},H^{\overline{\otimes}j})$
\end{center}

\noindent(a $2$-cocycle in the $3$-complex appears in Figure 3). The \emph{%
G-S }$n$-\emph{cohomology} \emph{of }$H$\emph{\ with coefficients in }$H$ is
given by

\begin{center}
$H_{GS}^{\ast}\left( H,H;n\right) =H_{\ast}\left\{ C_{GS}^{r}\left(
H,H;n\right) ;D\right\} $.
\end{center}

\noindent Note that a general $2$-cocycle $\alpha$ has a component of
tridegree $\left( 3-i-j,i,j\right) $ for each $i$ and $j$ in the range $%
2\leq i+j\leq n+1$. Thus $\alpha$ has $n\left( n+1\right) /2$ components and
a standard result in deformation theory tells us that the homogeneous
components of $\alpha$ determine an $\emph{infinitesimal}$ $\emph{deformation%
}$, i.e., the component $\omega_{1}^{j,i}$ in tridegree $\left(
3-i-j,i,j\right) $ defines the first order approximation $%
\omega_{0}^{j,i}+t\omega_{1}^{j,i}$ of the structure map $\omega_{t}^{j,i}$
in $H_{t}$.

For simplicity, consider the case $n=3.$ Each of the ten homogeneous
components of the deformation equation $D\left( \alpha\right) =0$ produces
the infinitesimal form of one structure relation (see below). In particular,
a deformation $H_{t}$ with structure maps $\{ \omega_{t}^{1,i}\} _{1\leq
i\leq3}$ is a simple $A\left( 3\right) $-algebra and a deformation $H_{t}$
with structure maps $\{ \omega_{t}^{j,1}\} _{1\leq j\leq3}$ is a simple $%
A\left( 3\right) $-coalgebra.

\vspace*{0.1in}

\begin{center}
\hspace*{0.5in}%
\includegraphics[
trim=0.000000in -0.222875in 2.231120in -0.222875in,
height=2.2in,
width=2.8in
]{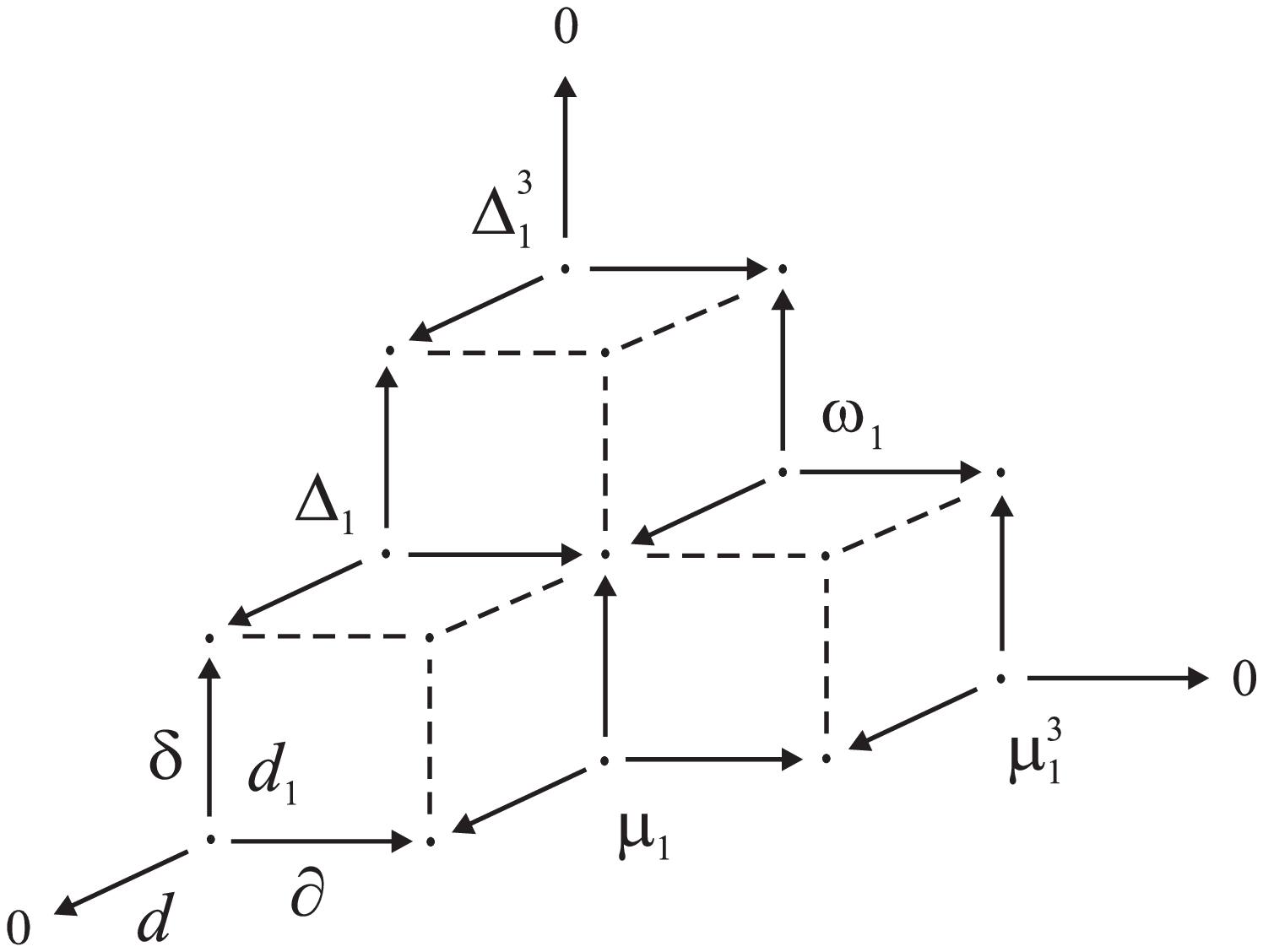}\\[0pt]
Figure 3. The 2-cocycle $d_{1}+\mu_{1}+\Delta_{1}+\mu_{1}^{3}+\omega_{1}+%
\Delta_{1}^{3}.$
\end{center}

\vspace*{0.1in}

For notational simplicity, let $\mu_{t}^{3}=\omega_{t}^{1,3},\,\ \omega
_{t}=\omega_{t}^{2,2}$ and $\Delta_{t}^{3}=\omega^{3,1},$ and consider a
deformation of $\left( H,d,\mu,\Delta\right) $ as a \textquotedblleft quasi-$%
A\left( 3\right) $-structure.\textquotedblright\ Then

\begin{itemize}
\item $d_{t}=d+td_{1}+t^{2}d_{2}+\cdots$

\item $\mu_{t}=\mu+t\mu_{1}+t^{2}\mu_{2}+\cdots$

\item $\Delta_{t}=\Delta+t\Delta_{1}+t^{2}\Delta_{2}+\cdots$

\item $\mu_{t}^{3}=t\mu_{1}^{3}+t^{2}\mu_{2}^{3}+\cdots$

\item $\omega_{t}=t\omega_{1}+t^{2}\omega_{2}+\cdots$

\item $\Delta_{t}^{3}=t\Delta_{1}^{3}+t^{2}\Delta_{2}^{3}+\cdots$
\end{itemize}

\noindent and $d_{1}+\mu_{1}+\Delta_{1}+\mu_{1}^{3}+\omega_{1}+%
\Delta_{1}^{3} $ is a total $2$-cocycle (see Figure 3). Equating
coefficients in $D\left(
d_{1}+\mu_{1}+\Delta_{1}+\mu_{1}^{3}+\omega_{1}+\Delta_{1}^{3}\right) =0$
gives

$
\begin{array}{llllll}
1. & d\left( d_{1}\right) =0 & \hspace*{0.5in} &  & 6. & \partial \left( \mu
_{1}^{3}\right) =0 \\ 
&  &  &  &  &  \\ 
2. & d\left( \mu _{1}\right) -\partial \left( d_{1}\right) =0 &  &  & 7. & 
\delta \left( \Delta _{1}^{3}\right) =0 \\ 
&  &  &  &  &  \\ 
3. & d\left( \Delta _{1}\right) +\delta \left( d_{1}\right) =0 &  &  & 8. & 
d\left( \omega _{1}\right) +\partial \left( \Delta _{1}\right) +\delta
\left( \mu _{1}\right) =0 \\ 
&  &  &  &  &  \\ 
4. & d\left( \mu _{1}^{3}\right) +\partial \left( \mu _{1}\right) =0 &  &  & 
9. & \partial \left( \Delta _{1}^{3}\right) +\delta \left( \omega
_{1}\right) =0 \\ 
&  &  &  &  &  \\ 
5. & d\left( \Delta _{1}^{3}\right) -\delta \left( \Delta _{1}\right) =0 & 
&  & 10. & \partial \left( \omega _{1}\right) -\delta \left( \mu
_{1}^{3}\right) =0.%
\end{array}%
$ 

Requiring $\left( H,d_{t},\mu_{t},\mu_{t}^{3}\right) $ and $\left(
H,d_{t},\Delta_{t},\Delta_{t}^{3}\right) $ to be simple $A\!\left( 3\right) $%
-(co)algebras tells us that relations (1) - (7) are linearizations of
Stasheff's strict $A\left( 4\right) $-(co)algebra relations, and relation
(8) is the linearization of the Hopf relation relaxed up to homotopy. Since $%
\mu_{t},$ $\omega_{t}$ and $\Delta_{t}$ have no terms of order zero,
relations (9) and (10) are the respective linearizations of new relations
(9) and (10) below. Thus we obtain the following structure relations in $%
H_{t}$:
\vspace{-0.1in}
\begin{center}
$%
\begin{array}{ll}
1. & d_{t}^{2}=0 \\ 
&  \\ 
2. & d_{t}\mu_{t}=\mu_{t}\left( d_{t}\otimes1+1\otimes d_{t}\right) \\ 
&  \\ 
3. & \Delta_{t}d_{t}=\left( d_{t}\otimes1+1\otimes d_{t}\right) \Delta_{t}
\\ 
&  \\ 
4. & d_{t}\mu_{t}^{3}+\mu_{t}^{3}\left( d_{t}\hspace{-0.02in}\otimes \hspace{%
-0.02in}1\hspace{-0.02in}\otimes\hspace{-0.02in}1+1\hspace {-0.02in}\otimes%
\hspace{-0.02in}d_{t}\hspace{-0.02in}\otimes\hspace {-0.02in}1+1\hspace{%
-0.02in}\otimes\hspace{-0.02in}1\hspace{-0.02in}\otimes\hspace{-0.02in}%
d_{t}\right) =\mu_{t}\left( 1\otimes\mu_{t}\right) -\mu_{t}\left(
\mu_{t}\otimes1\right) \\ 
&  \\ 
5. & \left( d_{t}\hspace{-0.02in}\otimes\hspace{-0.02in}1\hspace {-0.02in}%
\otimes\hspace{-0.02in}1+1\hspace{-0.02in}\otimes\hspace {-0.02in}d_{t}%
\hspace{-0.02in}\otimes\hspace{-0.02in}1+1\hspace{-0.02in}\otimes\hspace{%
-0.02in}1\hspace{-0.02in}\otimes\hspace{-0.02in}d_{t}\right)
\Delta_{t}^{3}+\Delta_{t}^{3}d_{t}=\left( \Delta_{t}\otimes1\right)
\Delta_{t}-\left( 1\otimes\Delta_{t}\right) \Delta_{t} \\ 
&  \\ 
6. & \mu_{t}^{3}\left( \mu_{t}\otimes1\otimes1-1\otimes\mu_{t}\otimes
1+1\otimes1\otimes\mu_{t}\right) =\mu_{t}\left( \mu_{t}^{3}\otimes
1+1\otimes\mu_{t}^{3}\right) \\ 
&  \\ 
7. & \left( \Delta_{t}\otimes1\otimes1-1\otimes\Delta_{t}\otimes
1+1\otimes1\otimes\Delta_{t}\right) \Delta_{t}^{3}=\left(
\Delta_{t}^{3}\otimes1+1\otimes\Delta_{t}^{3}\right) \Delta_{t} \\ 
&  \\ 
8. & \left( d_{t}\hspace{-0.02in}\otimes\hspace{-0.02in}1+1\hspace {-0.02in}%
\otimes\hspace{-0.02in}d_{t}\right) \omega_{t}+\omega_{t}\left( d_{t}\hspace{%
-0.02in}\otimes\hspace{-0.02in}1+1\hspace{-0.02in}\otimes \hspace{-0.02in}%
d_{t}\right) =\Delta_{t}\mu_{t}-\left( \mu_{t}\otimes\mu _{t}\right)
\sigma_{2,2}\left( \Delta_{t}\otimes\Delta_{t}\right) \\ 
&  \\
9. & \left( \mu_{t}\otimes\omega_{t}\right) \sigma_{2,2}\left( \Delta
_{t}\otimes\Delta_{t}\right) -\left( \Delta_{t}\otimes1-1\otimes\Delta
_{t}\right) \omega_{t}-\left( \omega_{t}\otimes\mu_{t}\right) \sigma
_{2,2}\left( \Delta_{t}\otimes\Delta_{t}\right) \\ 
&  \\ 
& \hspace*{0.3in}=\Delta_{t}^{3}\mu_{t}-\mu_{t}^{\otimes3}\sigma_{3,2}\left[
\left( \Delta_{t}\otimes1\right) \Delta_{t}\otimes\Delta_{t}^{3}+\left(
\Delta_{t}^{3}\otimes\left( 1\otimes\Delta_{t}\right) \Delta_{t}\right) %
\right] \\ 
&  \\
10. & \left( \mu_{t}\otimes\mu_{t}\right) \sigma_{2,2}\left( \Delta
_{t}\otimes\omega_{t}\right) -\omega_{t}\left( \mu_{t}\otimes1-1\otimes
\mu_{t}\right) -\left( \mu_{t}\otimes\mu_{t}\right) \sigma_{2,2}\left(
\omega_{t}\otimes\Delta_{t}\right) \\ 
&  \\ 
& \hspace*{0.3in}=\left[ \mu_{t}\left( \mu_{t}\otimes1\right) \otimes
\mu_{t}^{3}+\mu_{t}^{3}\otimes\mu_{t}\left( 1\otimes\mu_{t}\right) \right]
\sigma_{2,3}\Delta_{t}^{\otimes3}-\Delta_{t}\mu_{t}^{3}.  
\end{array}
$
\end{center}

\vspace*{0.1in}

\noindent By dropping the formal deformation parameter $t,$ we obtain the
structure relations in a \emph{quasi-simple} $A\left( 3\right) $-bialgebra.

The first \emph{non-operadic }example of an $A_{\infty}$-bialgebra appears
here as a quasi-simple $A\left( 3\right) $-bialgebra and involves a
non-trivial operation $\omega=\omega^{2,2}.$ The six additional relations
satisfied by $A_{\infty}$-bialgebras of this particular form will be
verified in the next section.

\begin{example}
\label{strict-A(4)}Let $H$ be the primitively generated bialgebra $%
\Lambda\left( x,y\right) $ with $\left\vert x\right\vert =1,$ $\left\vert
y\right\vert =2,$ trivial differential, and $\omega:H^{\otimes2}\rightarrow
H^{\otimes2}$ given by
\end{example}

\begin{center}
$\omega\left( a|b\right) =\left\{ 
\begin{array}{ll}
x|y+y|x, & a|b=y|y \\ 
x|x, & a|b\in\left\{ x|y,y|x\right\} \\ 
0, & \text{otherwise}.%
\end{array}
\right. $
\end{center}

\noindent Then $\left( \Delta\otimes1-1\otimes\Delta\right) \omega\left(
y|y\right) =\left( \Delta\otimes1-1\otimes\Delta\right) \left(
x|y+y|x\right) =1|x|y+1|y|x-x|y|1-y|x|1=\left( \mu\otimes\omega
-\omega\otimes\mu\right) \left( 1|1|y|y+y|1|1|y+1|y|y|1+y|y|1|1\right)
=\linebreak\left( \mu\otimes\omega-\omega\otimes\mu\right) \sigma
_{2,2}\left( \Delta\otimes\Delta\right) \left( y|y\right) $; similar
calculations show agreement on $x|y$ and $y|x$ and verifies relation (9). To
verify relation (10), note that $\omega\left( \mu\otimes1-1\otimes\mu\right) 
$ and $\left( \mu\otimes\mu\right) \sigma_{2,2}\left( \Delta\otimes
\omega-\omega\otimes\Delta\right) $ are supported on the subspace spanned by

\begin{center}
$B=\left\{ 1|y|y,y|y|1,1|x|y,x|y|1,1|y|x,y|x|1\right\} $,
\end{center}

\noindent and it is easy to check agreement on $B.$ Finally, note that $%
\left( H,\mu,\Delta,\omega\right) $ can be realized as the linear
deformation $\left. \left( H\left[ \left[ t\right] \right] ,\mu
,\Delta,t\omega\right) \right\vert _{t=1}.$

\section{$A_{\infty}$-Bialgebras in Perspective}

Although G-S deformation cohomology motivates the notion of an $A_{\infty}$%
-bialgebra, G-S deformations of DGBs are less constrained than $A_{\infty}$%
-bialgebras and fall short of the mark. To indicate of the extent of this
shortfall, let us identify those structure relations that fail to appear via
deformation cohomology but must be verified to assert that Example \ref%
{strict-A(4)} is an $A_{\infty}$-bialgebra.

As mentioned above, structure relations in a general $A_{\infty }$-bialgebra
arise from the homogeneous components of the equation $d_{\omega} \circledcirc d_{\omega}=0$. 
So to begin, let us construct the components of the biderivative $d_{\omega}$
that determine the structure relations in an $A_{\infty }$%
-bialgebra of the form $\left( H,d,\mu ,\Delta ,\omega ^{2,2}\right) $.
Given arbitrary maps $d=\omega ^{1,1},$ $\mu =\omega ^{1,2},$ $\Delta
=\omega ^{2,1},$ and $\omega ^{2,2}$ with $\omega ^{j,i}\in
Hom^{3-i-j}\left( H^{\otimes i},H^{\otimes j}\right) ,$ consider $\omega
=\sum \omega ^{j,i}.$ (Co)freely extend

\begin{itemize}
\item $d$ as a linear map $\left( H^{\otimes p}\right) ^{\otimes
q}\rightarrow\left( H^{\otimes p}\right) ^{\otimes q}$ for each $p,q\geq1,$

\item $d+\Delta$ as a derivation of $T^{a}H$,

\item $d+\mu$ as a coderivation of $T^{c}H$,

\item $\Delta +\omega ^{2,2}$ as a coalgebra map $T^{c}H\rightarrow
T^{c}\left( H^{\otimes 2}\right) $, and

\item $\mu +\omega ^{2,2}$ as an algebra map $T^{a}\left( H^{\otimes
2}\right) \rightarrow T^{a}H$.
\end{itemize}

\noindent Note that in this restricted setting, relation (10) in Definition %
\ref{A(4)-defn} reduces to

\begin{center}
$\left( \mu\otimes\mu\right) \circledcirc\left(
\Delta\otimes\omega^{2,2}-\omega^{2,2}\otimes\Delta\right)
=\omega^{2,2}\circledcirc\left( \mu \otimes1-1\otimes\mu\right) $.
\end{center}

\noindent Factors $\mu\otimes1$ and $1\otimes\mu$ are components of $%
\overline{d+\mu}$; factors $\Delta\otimes\omega^{2,2}$ and $\omega
^{2,2}\otimes\Delta$ are components of $\overline{\Delta+\omega^{2,2}}$; and
the factor $\mu\otimes\mu$ is a component of $\overline{\mu+\omega^{2,2}}$.

\vspace*{0.5in} \hspace*{1.6in}\setlength{\unitlength}{0.0002in} 
\begin{picture}(0,0)
\thicklines
\put(800,-3000){\line(1,0){5000}}
\put(800,-3000){\line(0,1){5000}}
\put(4800,-3000){\makebox(0,0){$\bullet$}}
\put(800,800){\makebox(0,0){$\bullet$}}
\put(4150,-2300){\vector(-1,1){2800}}
\put(-200,800){\makebox(0,0){$H^{\otimes j}$}}
\put(4800,-3600){\makebox(0,0){$H^{\otimes i}$}}
\put(3600,-200){\makebox(0,0){$\omega^{j,i}$}}
\end{picture}\vspace{0.8in}

\begin{center}
Figure 4. The initial map $\omega^{j,i}$.\vspace*{0.1in}
\end{center}

To picture this, identify the isomorphic modules $\left( H^{\otimes
p}\right) ^{\otimes q}\approx\left( H^{\otimes q}\right) ^{\otimes p}$ with
the point $\left( p,q\right) \in\mathbb{N}^{2}$ and picture the initial map $%
\omega^{j,i}:H^{\otimes i}\rightarrow H^{\otimes j}$ as a \textquotedblleft
transgressive\textquotedblright\ arrow\ from $\left( i,1\right) $ to $\left(
1,j\right) $ (see Figure 4). 

Components of the various (co)free extensions above are pictured as arrows
that initiate or terminate on the axes. For example, the vertical arrow $%
\Delta\otimes\Delta,$ the short left-leaning arrow $\Delta\otimes\omega
^{2,2}-\omega^{2,2}\otimes\Delta$ and the long left-leaning arrow $%
\omega^{2,2}\otimes\omega^{2,2}$ in Figure 5 represent components of $%
\overline{\Delta+\omega^{2,2}}$. 

\vspace*{0.1in}

\begin{center}
\includegraphics[
trim=-0.143316in -0.143051in -0.143316in -0.143050in,
height=2.7841in,
width=2.784in
]{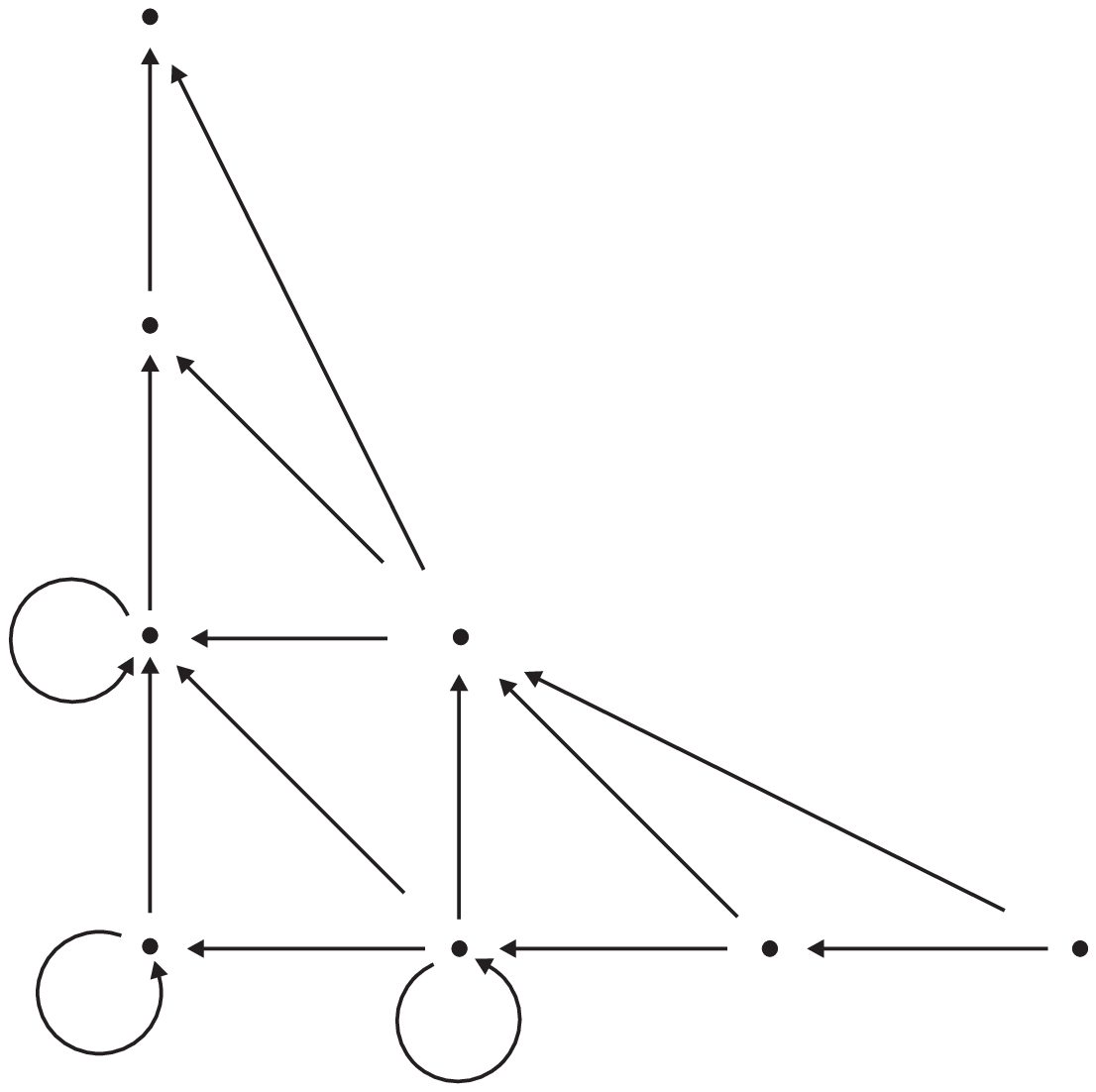}\\[0pt]
Figure 5. Components of $d_{\omega}$ when $\omega=d+\mu+\Delta+\omega^{2,2}.$
\end{center}

\vspace*{0.1in}

Since we are only interested in
transgressive quadratic $\circledcirc$-compositions, it is sufficient to
consider the components of $d_{\omega}$ pictured in Figure 5. Quadratic
compositions along the $x$-axis correspond to relations (1), (2), (4) and
(6) in Definition \ref{A(4)-defn}; those in the square with its diagonal
correspond to relation (8); those in the vertical parallelogram correspond
to relation (9); and those in the horizontal parallelogram correspond to
relation (10).

The following six additional relations are not detected by deformation
cohomology because the differentials only detect the interactions between $%
\omega$ and (deformations of) $d,$ $\mu,$ and $\Delta$ induced by the
underlying bi(co)module structure:$\vspace *{0.1in}$

\begin{enumerate}
\item[11.] $\left( \mu\otimes\omega-\omega\otimes\mu\right) \sigma
_{2,2}\left( \Delta\otimes\omega-\omega\otimes\Delta\right) =0;\vspace *{%
0.1in}$

\item[12.] $\left( \mu\otimes\mu\right) \sigma_{2,2}\left( \omega
\otimes\omega\right) =0;\vspace*{0.1in}$

\item[13.] $\left( \omega\otimes\omega\right) \sigma_{2,2}\left(
\Delta\otimes\Delta\right) =0;\vspace*{0.1in}$

\item[14.] $\left( \mu\otimes\omega-\omega\otimes\mu\right) \sigma
_{2,2}\left( \omega\otimes\omega\right) =0;\vspace*{0.1in}$

\item[15.] $\left( \omega\otimes\omega\right) \sigma_{2,2}\left(
\Delta\otimes\omega-\omega\otimes\Delta\right) =0;\vspace*{0.1in}$

\item[16.] $\left( \omega\otimes\omega\right) \sigma_{2,2}\left(
\omega\otimes\omega\right) =0.\vspace*{0.1in}$
\end{enumerate}

\begin{definition}
Let $H$ be a $\mathbf{k}$-module together with and a family of maps
$\left\{
d=\omega^{1,1},\mu=\omega^{1,2},\Delta=\omega^{2,1},\omega^{2,2}\right\}$, where
$\omega^{j,i}\in Hom^{3-i-j}\left( H^{\otimes i},H^{\otimes j}\right), $
and let $\omega=\sum\omega^{j,i}.$ Then $\left( H,d,\mu,\Delta,\omega
^{2,2}\right) $ is an $A_{\infty}$\textbf{-bialgebra} if $d_{\omega}\
\circledcirc\ d_{\omega}=0.$
\end{definition}

\begin{example}
\label{ex3}Continuing Example \ref{strict-A(4)}, verification of relations
(11) - (16) above is straightforward and follows from the fact that $%
\sigma_{2,2}\left( y|x|x|y\right) =-y|x|x|y.$ Thus $\left(
H,\mu,\Delta,\omega\right) $ is an $A_{\infty}$-bialgebra with non-operadic
structure.
\end{example}

Let $H$ be a graded module and let $\left\{ \omega ^{j,i}:H^{\otimes
i}\rightarrow H^{\otimes j}\right\} _{i,j\geq 1}$ be an arbitrary family of
maps. Given a diagonal $\Delta _{P}$ on the permutahedra and the notion of a 
$\Delta _{P}$-(co)derivation, one continues the procedure described above to
obtain the general biderivative defined in \cite{SU3}. And as above, the
general $A_{\infty }$-bialgebra structure relations are the homogeneous
components of $d_{\omega }\circledcirc d_{\omega }=0$.

For example, consider an $A_{\infty }$-bialgebra $(H,\mu ,\Delta ,\omega
^{j,i})$ with exactly one higher order operation $\omega ^{j,i}$, $i+j\geq 5.
$ When constructing $d_{\omega }$, we extend $\mu $ as a coderivation,
identify the components of this extension in $Hom\!\left( H^{\otimes
i}\!,H^{\otimes j}\!\right) $ with the vertices of the permutahedron $%
P_{i+j-2},$ and identify $\omega ^{j,i}$ with its top dimensional cell.
Since $\mu ,$ $\Delta $ and $\omega ^{j,i}$ are the only operations in $H$,
all compositions involving these operations have degree $0$ or $3-i-j$, and $%
k$-faces of $P_{i+j-2}$ in the range $0<k<i+j-3$ are identified with zero.
Thus the extension of $\omega ^{j,i}$ as a $\Delta _{P}$-coderivation only
involves the primitive terms of $\Delta \left( P_{i+j-2}\right) $, and the
components of this extension are terms in the expression $\delta \left(
\omega ^{j,i}\right) .$ Indeed, whenever $\omega ^{j,i}$ and its extension
are compatible with the underlying DGB structure, the relation $\delta
\left( \omega ^{j,i}\right) =0$ is satisfied. 

\begin{center}
\hspace{0.2in} 
\includegraphics[
trim=-0.143173in -0.143309in -0.143173in -0.143310in,
height=2.9983in,
width=4.0628in
]{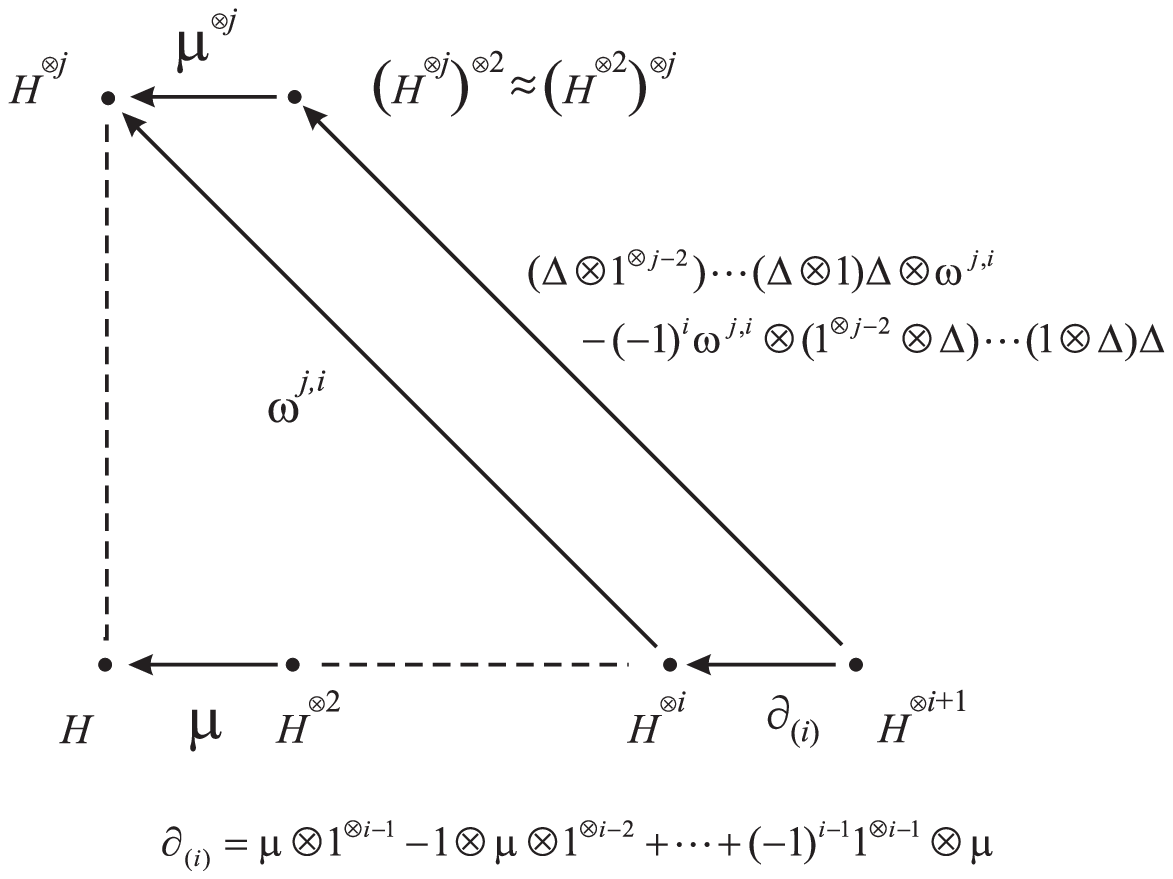}\\[0pt]
Figure 6. The structure relation $\partial\left( \omega^{j,i}\right) =0.$
\end{center}

\begin{center}
\hspace{0.3in} 
\includegraphics[
trim=-0.142990in -0.143244in -0.142991in -0.143244in,
height=3.286in,
width=3.4978in
]{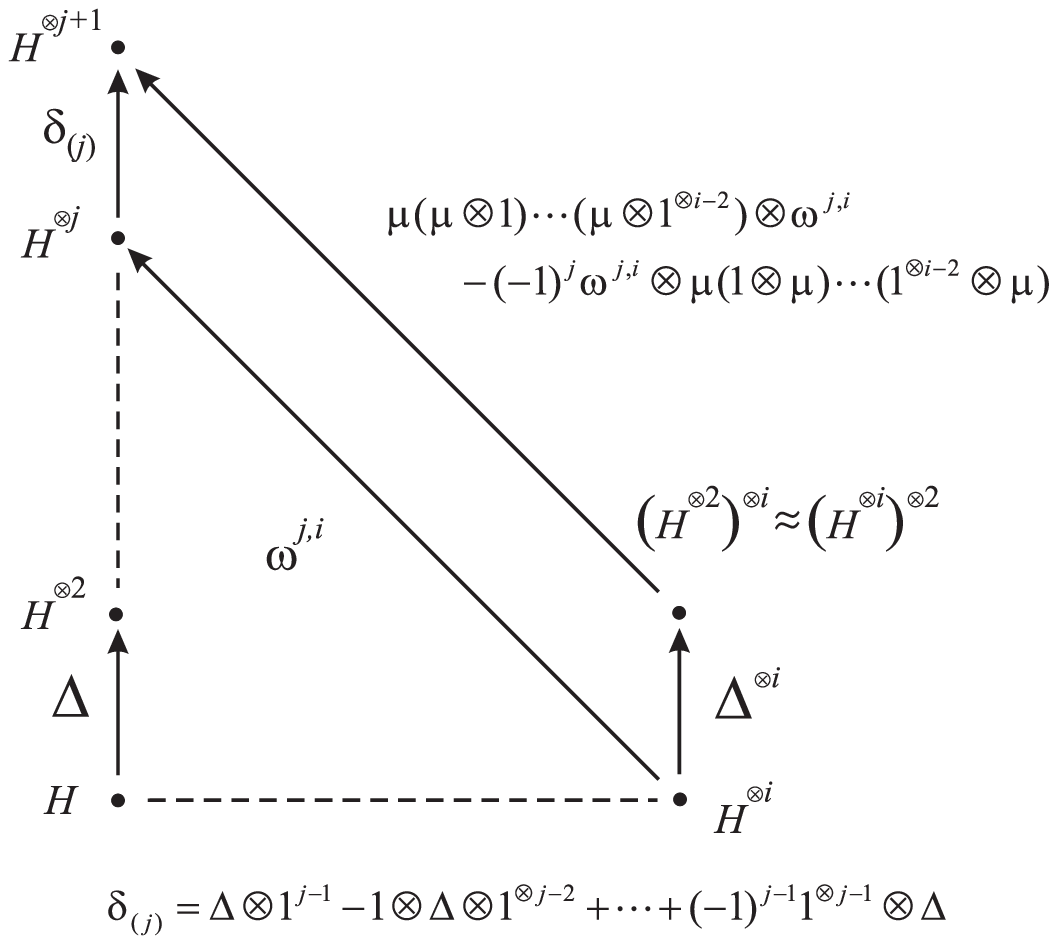}\\[0pt]
Figure 7. The structure relation $\delta\left( \omega^{j,i}\right) =0.$
\end{center}

Dually, we have $\partial
\left( \omega ^{j,i}\right) =0$ whenever $\omega ^{j,i}$ and its extension
as a $\Delta _{P}$-derivation are compatible with the underlying DGB
structure. These structure relations can be expressed as commutative
diagrams in the integer lattice $\mathbb{N}^{2}$ (see Figures 6 and 7).

\begin{definition}
\label{Hopf}Let $\ n\geq3.$ A \textbf{Hopf} $A\left( n\right) $\textbf{-algebra}
is a tuple $( H,d,\mu,\Delta,$\linebreak$\mu^{n}) $ with the following
properties:

\begin{enumerate}
\item[1.] $\left( H,d,\Delta\right) $ is a coassociative DGC;

\item[2.] $\left( H,d,\mu,\mu^{n}\right) $ is an $A\left( n\right) $%
-algebra; and

\item[3.] $\Delta\mu^{n}\!=\![\mu\left(
\mu\!\otimes\!1\right)\!\cdots\!\left( \mu\!\otimes\!1^{\otimes n-2}\right)
\otimes\mu^{n}+\mu^{n}\otimes\mu\left( 1\!\otimes\mu\right)\!\cdots\!\left(
1^{\otimes n-2}\!\otimes\mu\right) ]\sigma_{2,n}\Delta^{\otimes n}\!.$
\end{enumerate}

\noindent A \textbf{Hopf} $A_{\infty}$\textbf{-algebra} $\left( H,d,\mu
,\Delta,\mu^{n}\right) $ is a Hopf $A(n)$-algebra satisfying the
relation in offset (\ref{relation}) above. There are the completely dual
notions of a Hopf $A\left( n\right)$-coalgebra and a Hopf $%
A_{\infty}$-coalgebra.
\end{definition}

Hopf $A_{\infty}$-(co)algebras were defined by A. Berciano and this
author in \cite{Berciano}, but with a different choice of signs. $A_{\infty}$-bialgebras with operations
exclusively of the forms $\omega^{j,1}$ and $\omega^{1,i},$ called \emph{%
special} $A_{\infty}$\emph{-bialgebras, }were considered by this author in 
\cite{Umble2}. 

Hopf $A\left( n\right) $-algebras are especially
interesting because their structure relations are controlled by G-S
deformation cohomology. In fact, if $n\geq3$ and $H_{t}=\left( H\left[ \left[ t%
\right] \right] ,d_{t},\mu_{t},\Delta_{t},\mu_{t}^{n}\right) $ is a
deformation, then $\mu_{t}^{n}=t\mu_{1}^{n}+t^{2}\mu_{2}^{n}+\cdots$ has no
term of order zero. Consequently, if $D\left( \mu_{1}^{n}\right) =0,$ then $%
t\mu_{1}^{n}$ automatically satisfies the required structure relations
and $\left( H\left[ \left[ t\right] %
\right] ,d,\mu,\Delta,t\mu_{1}^{n}\right) $ is a \emph{linear} deformation
of $H_{0}$ as a Hopf $A\left( n\right) $-algebra. Thus we have proved:

\begin{theorem}
\label{main} \textit{If} $\left( H,d,\mu ,\Delta\right) $\textit{\ is a DGB
and\ }$\mu_{1}^{n}\in Hom^{2-n}\left( H^{\otimes n},H\right),$\textit{\ }$%
n\geq3,$\textit{\ is a 2-cocycle, then} $\left( H\left[ \left[ t\right] %
\right] ,d,\mu,\Delta ,t\mu_{1}^{n}\right) $ \textit{is a linear deformation
of }$H_{0}$ \textit{as a Hopf }$A\left( n\right) $-\textit{algebra.}
\end{theorem}

I am grateful to Samson Saneblidze and Andrey Lazarev for their helpful
suggestions on early drafts of this paper, and to the referee, the editors,
and Jim Stasheff for their assistance with the final draft. I wish to thank
Murray and Jim for their encouragement and support of this project over the
years and I wish them both much happiness and continued success. \bigskip %



\printindex

\end{document}